\documentclass{conm-p-l}

\usepackage{amssymb}
\usepackage[usenames, dvipsnames]{color}
\usepackage{tikz}
\usetikzlibrary{matrix,arrows}
\usepackage{blkarray}

\usepackage{}

\newcommand{\bd}{\mathbf{d}}
\DeclareMathOperator{\Mat}{Mat}
\DeclareMathOperator{\Hom}{Hom}
\DeclareMathOperator{\rank}{rank}

\newcommand{\rep}{\mathtt{rep}}
\newcommand{\cO}{\mathcal{O}}
\newcommand{\cOc}[1]{\overline{\mathcal{O}_{#1}}}
\newcommand{\bw}{\mathbf{w}}
\newcommand{\kk}{\Bbbk}

\newtheorem{theorem}{Theorem}[section]

\newtheorem{proposition}[theorem]{Proposition}

\theoremstyle{definition}

\theoremstyle{remark}
\newtheorem{remark}[theorem]{Remark}
\newtheorem{problem}[theorem]{Problem}

\numberwithin{equation}{section}

\begin{document}

\title{K-polynomials of type A quiver orbit closures and lacing diagrams}


\author{Ryan Kinser}
\address{University of Iowa, Department of Mathematics, Iowa City, IA, USA}
\email{ryan-kinser@uiowa.edu}
\thanks{}



\date{}

\begin{abstract}
This article contains an overview of the author's joint work with Allen Knutson and Jenna Rajchgot on $K$-polynomials of orbit closures for type $A$ quivers.  It is written to an audience interested in interactions between representations of algebras, algebraic geometry, and commutative algebra.  A few open problems resulting from the work are also explained.
\end{abstract}

\maketitle


\section{Background and context}
We denote a \emph{quiver} by $Q=(Q_0, Q_1, t, h)$, where $Q_0$ is the vertex set, $Q_1$ the arrow set, and $t, h\colon Q_1 \to Q_0$ give the tail and head of an arrow $ta \xrightarrow{a} ha$.  We fix a field $\kk$, and will proceed with the assumption that $\kk$ is algebraically closed in order to simplify the geometric language, although this is completely inessential to our main result since all schemes involved are defined over $\mathbb{Z}$.

Given a quiver $Q$ and \emph{dimension vector} $\bd \colon Q_0 \to \mathbb{Z}_{\geq 0}$, we study the \emph{representation variety}
\begin{equation}
\rep_Q(\bd) = \prod_{a \in Q_1} \Mat(\bd(ta), \bd(ha)),
\end{equation}
where $\Mat(m,n)$ denotes the variety of matrices with $m$ rows, $n$ columns, and entries in $\kk$.  
We consider the right action of the \emph{base change group}
\begin{equation}
GL(\bd) = \prod_{z \in Q_0} GL(\bd(z))
\end{equation}
on $\rep_Q(\bd)$ given by
\begin{equation}\label{eq:gM}
M\cdot g = (g_{ta}^{-1}M_a g_{ha})_{a\in Q_1},
\end{equation}
where $g = (g_z)_{z \in Q_0} \in GL(\bd)$ and $M = (M_a)_{a \in Q_1} \in \rep_Q(\bd)$.
A \emph{representation} of $Q$ is a collection of (finite-dimensional) $\kk$-vector spaces $(V_z)_{z \in Q_0}$ assigned to the vertices of $Q$, along with a collection of $\kk$-linear maps $(\varphi_a \colon V_{ta} \to V_{ha})_{a \in Q_1}$ assigned to the arrows.  Thus, the points of $\rep_Q(\bd)$ are in bijection with representations of $Q$ along with a fixed basis at each vertex.  For algebraic context, we mention that there is a natural definition of a \emph{morphism} between two representations which yields a category $\rep(Q)$ of all representations of $Q$. This category is abelian and in fact equivalent to the category of right modules over the \emph{path algebra} $\kk Q$.  Then, at least when $Q$ has no oriented cycles so that $\kk Q$ is finite dimensional, the dimension vector of a representation is equivalent to its class in the Grothendieck group of $\rep(Q)$.  We refer the interested reader to standard references such as \cite{Schiffler:2014aa,assemetal,ARSartinalgebras} for further details.

Then simply from the definitions, orbits in $\rep_Q(\bd)$ under $GL(\bd)$ are in bijection with isomorphism classes of representations of dimension vector $\bd$; for a representation $M$ of $Q$, we denote by $\cO_M$ the orbit of $M$ in $\rep_Q(\bd)$, and by $\cOc{M}$ the closure of this orbit.
Orbit closures in $\rep_Q(\bd)$ have remarkable connections with the representation theory of $Q$ and related objects. 
Here we highlight a few connections, and refer the reader to surveys such as \cite{Bongartzsurvey,Zwarasurvey,HZsurvey} for detailed treatments of the connections to representation theory.

\subsubsection*{Commutative algebra} 
Orbit closures in $\rep_Q(\bd)$ come with a natural set of polynomials vanishing on them obtained from projective resolutions of the indecomposable representations of $Q$.  These polynomials are minors of certain matrices whose entries are the natural coordinate functions on $\rep_Q(\bd)$, possibly repeated, and 0s (see \cite[\S4]{MR3008913} or \cite[\S3]{KR}).  So ideals generated by these minors can be seen as generalizations of determinantal ideals.  From this perspective it is then natural to ask when these ideals are primary, prime, normal, Cohen-Macaulay, etc.  There are some surprisingly general results, such as the fact that the ideals obtained in this way are always primary \cite{MR1402728} when $Q$ is a quiver of Dynkin type $A, D$, or $E$ (generalized to all representation-finite algebras by Zwara \cite{MR1476404}).  

\subsubsection*{Lie theory} Each Dynkin quiver $Q$ determines a finite-dimensional, simple complex Lie algebra $\mathfrak{g}_Q$ (which is actually independent of the orientation of $Q$) and thus a universal enveloping algebra $U(\mathfrak{g}_Q)$.  Ringel's work constructing the upper half $U(\mathfrak{n})$ of this algebra as a Hall algebra \cite{Rhallalgebras} was geometrically realized by Lusztig \cite{MR1035415,Lusztig:1991yq} as a convolution algebra of constructible functions on $\rep_Q(\bd)$ which are constant on the orbit closures (there was also an unpublished manuscript by Schofield on the subject around this time).  More recently, Geiss, Leclerc, and Schr\"oer have made strides towards generalizing some of the above mentioned work to arbitrary symmetrizable Kac-Moody Lie algebras (not necessarily simply laced) \cite{MR3555157}.  Their work uses convolution algebras of constructible functions which are constant on orbit closures in representation schemes of certain Iwanaga-Gorenstien algebras of dimension 1.

\subsubsection*{Representations of algebras} 
We will just mention here two kinds of results relating the representation theory of a quiver to the geometry of orbit closures in its representation varieties, and refer the interested reader to the surveys cited above for more.
The \emph{degeneration order} on representations of $Q$ (of the same dimension vector) is defined by $M \leq_{deg} N$ if and only if $\cOc{M} \supseteq \cOc{N}$.  
It turns out that this order is closely connected to algebraic properties of $M$ and $N$.  The nicest results are for $Q$ of Dynkin or extended Dynkin type, where for example $M \leq_{deg} N$ if and only if $\dim \Hom_Q(M, X) \leq \dim \Hom_Q(N, X)$ for all indecomposable representations $X$ \cite{MR1341664}.  The latter condition is called the \emph{Hom order} and typically denoted simply by $M \leq N$; this was further proven to be equivalent to a related \emph{Ext order} for extended Dynkin quivers (and more generally, tame concealed algebras) by Zwara \cite{MR1610587}.

Another remarkable characterization due to Zwara \cite{MR1757882}, building on work of Riedtmann \cite{MR868301}, is that $M \leq_{deg} N$ exactly when there exists another representation $Z$ (of unknown dimension) and an exact sequence of the form
\[
0 \to N \to M \oplus Z \to Z \to 0,
\]
and that in this case the degeneration can be realized by a rational curve $\mathbb{A}^1 \to \cOc{M}$.  

\subsubsection*{Algebraic geometry} The work described in this article is most directly inspired by the literature on degeneracy loci.  To start with the simplest case, given a nonsingular algebraic variety $X$ and a map between vector bundles $\phi\colon E \to F$ on $X$, let $\phi_x \colon E_x \to F_x$ denote the induced map on fibers over $x \in X$.
Then for $r \in \mathbb{Z}_{\geq 0}$, we may consider the \emph{degenerarcy locus} $\Omega_r = \{ x \in X \mid \rank \phi_x \leq r \}$, which is a closed subvariety of $X$ since it is defined by the vanishing of minors in local coordinates around each point.  Then it turns out that, when $\phi$ is sufficiently general, the fundamental class of $\Omega_r$ in the Chow or cohomology ring of $X$ has a universal expression as a Schur function evaluated at the Chern roots of $E$ and $F$ (the Giambelli-Thom-Porteous determinantal formula).  A brief history tracing this formula back to its geometric and algebraic roots in the 1800s can be found in the Introduction of Fulton's paper \cite{Fulton92}; this paper established many ideas essential to our work.

The connection with quivers originated with Buch and Fulton \cite{MR1671215,BFchernclass}, generalizing to sequences of vector bundle maps.  In quiver language, viewing $\phi \colon E \to F$ in the setup above as a ``representation of the quiver $A_2$ by vector bundles on $X$'', their work replaces $A_2$ by an arbitrary \emph{equioriented} type $A$ quiver (all arrows pointing in the same direction).  It is natural from there to seek generalizations to other quivers, although the strongest results are to be expected for Dynkin quivers.  Buch elevated the formulas to the level of $K$-theory in \cite{MR1932326}; 
see also \cite{MR2114821,FRdegenlocithom,BFR, KMS,RimanyiCOHA,MR3286674,MR3239295} for other important contributions.
Much of the state of the art can be found in Buch's article \cite{MR2492443}, which focuses on Dynkin quivers.

\smallskip

The aim of this article is to give a somewhat self-contained overview of the ``$K$-theoretic component formula'' proven in the author's joint work with Allen Knutson and Jenna Rajchgot \cite{KKR}.
Readers interested in the interface between representation theory of algebras and algebraic geometry may find this article to be a more accessible introduction to the result, as it contains a little more commentary and references to the representation theory side than the original paper.  It also contains a single running example illustrating most of the key ideas.  In particular, we highlight the role of lacing diagrams in the equivariant geometry of orbit closures in $\rep_Q(\bd)$, since these should be intuitive to anyone familiar with representations of type $A$ quivers.  We shall also pose some open questions with the hope of motivating further work on the topic.

\section{Lacing diagrams}\label{sec:lacing}
Assume $Q$ is a quiver of Dynkin type $A$ for this section.  Lacing diagrams were introduced by Abeasis and del Fra in \cite{AdF} as a tool to combinatorially characterize the degeneration order on orbits in $\rep_Q(\bd)$ (which they call the ``geometrical ordering'').  Knutson, Miller, and Shimozono introduced a refinement of Abeasis and del Fra's diagrams \cite{KMS} for equioriented type $A$ quivers, realizing that allowing the laces to cross enables us to keep track of subtle combinatorial information relevant to the equivariant geometry of orbit closures in $\rep_Q(\bd)$.
Buch and Rim\'anyi utilized lacing diagrams to study equivariant geometry of $\rep_Q(\bd)$ for arbitrarily oriented $Q$ of type $A$ in \cite{MR2306279}.

A \emph{lacing diagram} of dimension vector $\bd$ for $Q$ consists of:
\begin{enumerate}
\item for each vertex $z \in Q_0$, a column of $\bd(z)$ dots;
\item for each arrow $a \in Q_1$, a set of arrows from dots in column $ta$ to dots in column $ha$, such that no dot has more than one incoming or outgoing arrow.
\end{enumerate}

We will develop a running example throughout this article using the following quiver.
\begin{equation}
Q=\quad
\vcenter{\hbox{\begin{tikzpicture}[point/.style={shape=circle,fill=black,scale=.5pt,outer sep=3pt},>=latex]
   \node[outer sep=-2pt] (1) at (2,1) {${4}$};
  \node[outer sep=-2pt] (2) at (0,1) {${2}$};
   \node[outer sep=-2pt] (3) at (1,0) {${3}$};
   \node[outer sep=-2pt] (4) at (-1,0) {${1}$};
  \path[->]
	(2) edge node[auto] {${\alpha}$} (4)
	(2) edge node[auto] {${\beta}$} (3)
  	(1) edge node[auto] {${\gamma}$} (3);
   \end{tikzpicture}}}
\end{equation}
Three lacing diagrams of dimension vector $\bd = (\bd(1), \bd(2), \bd(3), \bd(4)) = ( 1,3,2,1)$ for $Q$ are seen in Figure \ref{fig:laces}.

\begin{figure}[b]
\begin{tikzpicture}[point/.style={shape=circle,fill=black,scale=.5pt,outer sep=3pt},epoint/.style={shape=rectangle,fill=red,scale=.5pt,outer sep=3pt},>=latex] 
\node[point] (1a) at (0,0) {};
\node[point] (2a) at (1,0) {};
\node[point] (2b) at (1,1) {};
\node[point] (2c) at (1,2) {};
\node[point] (3b) at (2,1) {};
\node[point] (3c) at (2,2) {};
\node[point] (4b) at (3,1) {};
  
\path[->,thick]
   (2c) edge (1a)
   (2b) edge (3b)
   (2c) edge (3c)
   (4b) edge (3b);  
  \end{tikzpicture}
\qquad
\begin{tikzpicture}[point/.style={shape=circle,fill=black,scale=.5pt,outer sep=3pt},epoint/.style={shape=rectangle,fill=red,scale=.5pt,outer sep=3pt},>=latex] 
\node[point] (1a) at (0,0) {};
\node[point] (2a) at (1,0) {};
\node[point] (2b) at (1,1) {};
\node[point] (2c) at (1,2) {};
\node[point] (3b) at (2,1) {};
\node[point] (3c) at (2,2) {};
\node[point] (4b) at (3,1) {};
  
\path[->,thick]
   (2b) edge (1a)
   (2b) edge (3c)
   (2c) edge (3b)
   (4b) edge (3b);  
  \end{tikzpicture}
\qquad
\begin{tikzpicture}[point/.style={shape=circle,fill=black,scale=.5pt,outer sep=3pt},epoint/.style={shape=rectangle,fill=red,scale=.5pt,outer sep=3pt},>=latex] 
\node[point] (1a) at (0,0) {};
\node[point] (2a) at (1,0) {};
\node[point] (2b) at (1,1) {};
\node[point] (2c) at (1,2) {};
\node[point] (3b) at (2,1) {};
\node[point] (3c) at (2,2) {};
\node[point] (4b) at (3,1) {};
  
\path[->,thick]
   (2b) edge (1a)
   (2b) edge (3b)
   (2a) edge (3c)
   (4b) edge (3c);  
  \end{tikzpicture}
    \caption{Three lacing diagrams for $Q$.}
\label{fig:laces}
\end{figure}

Equivalently, one may simply define a lacing diagram as a sequence $\bw = (w_a)_{a \in Q_1}$ where each $w_a$ is a $\bd(ta) \times 
\bd(ha)$ \emph{partial permutation matrix}, meaning its entries are all 0 or 1, with at most one 1 in each row and each column.  
Note that this definition clearly identifies a lacing diagram $\bw$ as a specific point of $\rep_Q(\bd)$.
For details of our conventions, which may differ from other authors for technical reasons, see \cite[\S2.8]{KKR}.  The matrix representation of the leftmost lacing diagram in Figure \ref{fig:laces} is
\begin{equation}
\bw = (w_\alpha, w_\beta, w_\gamma) = 
\left(
\begin{bmatrix}
1 \\
0 \\
0
\end{bmatrix},
\begin{bmatrix}
1 & 0\\
0 & 1 \\
0 & 0
\end{bmatrix},
\begin{bmatrix}
0 & 1
\end{bmatrix}
\right).
\end{equation}

Note that the diagrammatic presentation of a lacing diagram makes it very easy to recognize the orbit of the point in $\rep$: this is equivalent to knowing its direct sum decomposition into indecomposables, which are just the individual laces.  For example, one clearly sees that all three lacing diagrams above lie in the same orbit because they have the same number of laces connecting any two columns.  They are also more convenient for performing certain combinatorial manipulations described below.  On the other hand, the partial permutation matrix viewpoint provides a more natural interface with combinatorial commutative algebra.

\begin{remark}
One possible way to generalize lacing diagrams to arbitrary quivers would be to consider \emph{tree modules}, or more precisely, the \emph{coefficient quivers} of tree modules \cite[\S2]{Ringel:1998gf}.  Roughly speaking, a tree module is an indecomposable quiver representation whose structure can be completely encoded by another quiver whose underlying graph is a tree.
Ringel has shown that any indecomposable quiver representation with no self extensions is a tree module, so in particular any indecomposable representation of a Dynkin quiver is a tree module, and therefore an arbitrary representation of a Dynkin quiver can be presented as a disjoint union of trees.  For type $A$ quivers, these specialize exactly to lacing diagrams.
Tree modules have been studied extensively by Ringel, Weist and others, see for example \cite{MR1090218,kinserrootedtrees,MR2578596,MR2975149,MR3102957,MR2988715,MR3183888}.
\end{remark}

Although the example shows that an orbit may be represented by many different lacing diagrams, we will see now that there are certain distinguished lacing diagrams.
Each lacing diagram has hidden ``virtual'' laces which must be revealed in order to read off the combinatorial information encoded in the diagram.   To do this, we need to fix another convention: let $\{1, 2, \dotsc, n\}$ denote the vertex set of $Q$, such that there is an arrow between $i$ and $j$ (in either direction) if and only $|i - j| = 1$.  Arrows of the form $i \xrightarrow{a} i+1$ will be called rightward arrows, and the others leftward arrows.  We complete a lace diagram $\bw = (w_a)$ to its \emph{extended lacing diagram} $(c_a(w_a))$ by extending each partial permutation matrix to a full permutation matrix according to the following convention: if $a$ is rightward pointing, then $c_a(w_a)$ is the unique permutation of minimal size and Coxeter length such that $w_a$ lies in the northwest corner of $c_a(w_a)$, and for $a$ leftward pointing, $c_a(w_a)$ is the same but containing $w_a$ in the southeast corner instead.  The length $|\bw|$ of $\bw$ is defined as
\[
|\bw| = \sum_{a \in Q_1} \ell(c_a(w_a)),
\]
where the function $\ell$ gives the Coxeter length of a permutation with respect to the standard generating set of adjacent transpositions.  A \emph{minimal lacing diagram} is one whose length is minimal among all lacing diagrams in $\cO_{\bw}$.  The completion is encoded in a visual diagram in a natural way by adding ``virtual'' dots and laces (red squares and dashed lines in our examples).

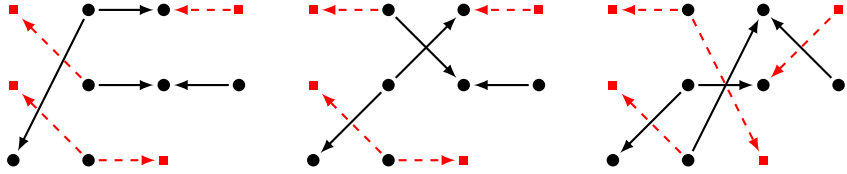
\begin{figure}
\begin{tikzpicture}[point/.style={shape=circle,fill=black,scale=.5pt,outer sep=3pt},epoint/.style={shape=rectangle,fill=red,scale=.5pt,outer sep=3pt},>=latex] 
\node[point] (1a) at (0,0) {};
\node[epoint] (1b) at (0,1) {};
\node[epoint] (1c) at (0,2) {};
\node[point] (2a) at (1,0) {};
\node[point] (2b) at (1,1) {};
\node[point] (2c) at (1,2) {};
\node[epoint] (3a) at (2,0) {};
\node[point] (3b) at (2,1) {};
\node[point] (3c) at (2,2) {};
\node[point] (4b) at (3,1) {};
\node[epoint] (4c) at (3,2) {};
  
\path[->,thick]
   (2a) edge[red,dashed] (1b)
   (2b) edge[red,dashed] (1c)
   (2c) edge (1a)
   (2a) edge[red,dashed] (3a)
   (2b) edge (3b)
   (2c) edge (3c)
   (4b) edge (3b)
   (4c) edge[red,dashed] (3c);
  \end{tikzpicture}
\qquad
\begin{tikzpicture}[point/.style={shape=circle,fill=black,scale=.5pt,outer sep=3pt},epoint/.style={shape=rectangle,fill=red,scale=.5pt,outer sep=3pt},>=latex] 
\node[point] (1a) at (0,0) {};
\node[epoint] (1b) at (0,1) {};
\node[epoint] (1c) at (0,2) {};
\node[point] (2a) at (1,0) {};
\node[point] (2b) at (1,1) {};
\node[point] (2c) at (1,2) {};
\node[epoint] (3a) at (2,0) {};
\node[point] (3b) at (2,1) {};
\node[point] (3c) at (2,2) {};
\node[point] (4b) at (3,1) {};
\node[epoint] (4c) at (3,2) {};
  
\path[->,thick]
   (2a) edge[red,dashed] (1b)
   (2c) edge[red,dashed] (1c)
   (2b) edge (1a)
   (2b) edge (3c)
   (2c) edge (3b)
   (2a) edge[red,dashed] (3a)
   (4b) edge (3b)
   (4c) edge[red,dashed] (3c);

  \end{tikzpicture}
\qquad
\begin{tikzpicture}[point/.style={shape=circle,fill=black,scale=.5pt,outer sep=3pt},epoint/.style={shape=rectangle,fill=red,scale=.5pt,outer sep=3pt},>=latex] 
\node[point] (1a) at (0,0) {};
\node[epoint] (1b) at (0,1) {};
\node[epoint] (1c) at (0,2) {};
\node[point] (2a) at (1,0) {};
\node[point] (2b) at (1,1) {};
\node[point] (2c) at (1,2) {};
\node[epoint] (3a) at (2,0) {};
\node[point] (3b) at (2,1) {};
\node[point] (3c) at (2,2) {};
\node[point] (4b) at (3,1) {};
\node[epoint] (4c) at (3,2) {};
  
\path[->,thick]
   (2a) edge[red,dashed] (1b)
   (2b) edge (1a)
   (2c) edge[red,dashed] (1c)
   (2b) edge (3b)
   (2a) edge (3c)
   (2c) edge[red,dashed] (3a)
   (4c) edge[red,dashed] (3b)
   (4b) edge (3c);  
  \end{tikzpicture}
    \caption{Three completed lacing diagrams for $Q$.}
\label{fig:laces2}
\end{figure}

The completions of the lacing diagrams in Figure \ref{fig:laces} are shown in Figure \ref{fig:laces2}.
The first two are minimal, with 2 crossings each, but the last is not since it lies in the same orbit but the completion has 5 crossings. 

The following proposition gives a first taste of the connections between the combinatorics of lacing diagrams, equivariant geometry of $\rep_Q(\bd)$, and representation theory of $Q$. It follows from \cite[Cor. 2]{MR2306279} and the Artin-Voigt formula (see, for example, \cite[Lemma 2.3]{Ringel:1980ly}).

\begin{proposition}
Let $\bw$ be a minimal lacing diagram in $\rep_Q(\bd)$.  Then $|\bw|$ is equal to both  the codimension of $\cOc{\bw}$ in $\rep_Q(\bd)$, and $\dim_\kk {\rm Ext}_Q^1(\bw, \bw)$.
\end{proposition}

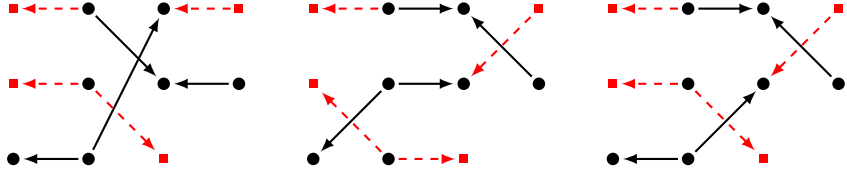
\begin{figure}[h]
\begin{tikzpicture}[point/.style={shape=circle,fill=black,scale=.5pt,outer sep=3pt},epoint/.style={shape=rectangle,fill=red,scale=.5pt,outer sep=3pt},>=latex] 
\node[point] (1a) at (0,0) {};
\node[epoint] (1b) at (0,1) {};
\node[epoint] (1c) at (0,2) {};
\node[point] (2a) at (1,0) {};
\node[point] (2b) at (1,1) {};
\node[point] (2c) at (1,2) {};
\node[epoint] (3a) at (2,0) {};
\node[point] (3b) at (2,1) {};
\node[point] (3c) at (2,2) {};
\node[point] (4b) at (3,1) {};
\node[epoint] (4c) at (3,2) {};
  
\path[->,thick]
   (2a) edge (1a)
   (2b) edge[red,dashed] (1b)
   (2c) edge[red,dashed] (1c)
   (2a) edge (3c)
   (2b) edge[red,dashed] (3a)
   (2c) edge (3b)
   (4b) edge (3b)
   (4c) edge[red,dashed] (3c);
  \end{tikzpicture}
\qquad
\begin{tikzpicture}[point/.style={shape=circle,fill=black,scale=.5pt,outer sep=3pt},epoint/.style={shape=rectangle,fill=red,scale=.5pt,outer sep=3pt},>=latex] 
\node[point] (1a) at (0,0) {};
\node[epoint] (1b) at (0,1) {};
\node[epoint] (1c) at (0,2) {};
\node[point] (2a) at (1,0) {};
\node[point] (2b) at (1,1) {};
\node[point] (2c) at (1,2) {};
\node[epoint] (3a) at (2,0) {};
\node[point] (3b) at (2,1) {};
\node[point] (3c) at (2,2) {};
\node[point] (4b) at (3,1) {};
\node[epoint] (4c) at (3,2) {};
  
\path[->,thick]
   (2a) edge[red,dashed] (1b)
   (2b) edge (1a)
   (2c) edge[red,dashed] (1c)
   (2a) edge[red,dashed] (3a)
   (2b) edge (3b)
   (2c) edge (3c)
   (4b) edge (3c)
   (4c) edge[red,dashed] (3b);

  \end{tikzpicture}
\qquad
\begin{tikzpicture}[point/.style={shape=circle,fill=black,scale=.5pt,outer sep=3pt},epoint/.style={shape=rectangle,fill=red,scale=.5pt,outer sep=3pt},>=latex] 
\node[point] (1a) at (0,0) {};
\node[epoint] (1b) at (0,1) {};
\node[epoint] (1c) at (0,2) {};
\node[point] (2a) at (1,0) {};
\node[point] (2b) at (1,1) {};
\node[point] (2c) at (1,2) {};
\node[epoint] (3a) at (2,0) {};
\node[point] (3b) at (2,1) {};
\node[point] (3c) at (2,2) {};
\node[point] (4b) at (3,1) {};
\node[epoint] (4c) at (3,2) {};
  
\path[->,thick]
   (2a) edge (1a)
   (2b) edge[red,dashed] (1b)
   (2c) edge[red,dashed] (1c)
   (2a) edge (3b)
   (2b) edge[red,dashed] (3a)
   (2c) edge (3c)
   (4b) edge (3c)
   (4c) edge[red,dashed] (3b);
  \end{tikzpicture}
    \caption{The remaining minimal lacing diagrams for $\cO_{\bw}$}
\label{fig:laces3}
\end{figure}

All minimal lacing diagrams in a given orbit are related by a sequence of moves of the following form, where both dots in the middle column and at least one dot in each outer column must not be virtual.
\begin{equation}
\vcenter{\hbox{  
\begin{tikzpicture}[point/.style={shape=circle,fill=black,scale=.5pt,outer sep=3pt},>=latex]
\node[point] (1a) at (0,0) {};
\node[point] (1b) at (0,1) {};
\node[point] (2a) at (1,0) {};
\node[point] (2b) at (1,1) {};
\node[point] (3a) at (2,0) {};
\node[point] (3b) at (2,1) {};

\path[-,thick]
   (1a) edge (2b)
   (1b) edge (2a)
   (2a) edge (3a)
   (2b) edge (3b);
\end{tikzpicture}}}
\longleftrightarrow
\vcenter{\hbox{  
\begin{tikzpicture}[point/.style={shape=circle,fill=black,scale=.5pt,outer sep=3pt},>=latex]
\node[point] (1a) at (0,0) {};
\node[point] (1b) at (0,1) {};
\node[point] (2a) at (1,0) {};
\node[point] (2b) at (1,1) {};
\node[point] (3a) at (2,0) {};
\node[point] (3b) at (2,1) {};

\path[-,thick]
   (1a) edge (2a)
   (1b) edge (2b)
   (2a) edge (3b)
   (2b) edge (3a);
\end{tikzpicture}}}
\end{equation}
The remaining minimal lacing diagrams for the particular orbit in the running example are in Figure \ref{fig:laces3}, so there are 5 total.

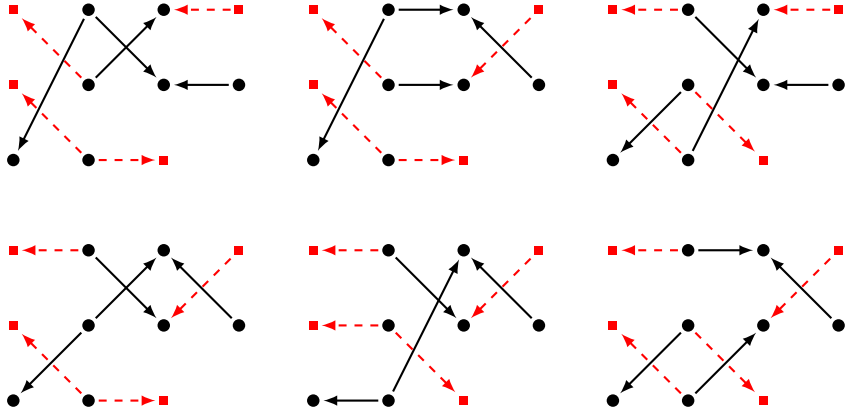
\begin{figure}[h]
\begin{tikzpicture}[point/.style={shape=circle,fill=black,scale=.5pt,outer sep=3pt},epoint/.style={shape=rectangle,fill=red,scale=.5pt,outer sep=3pt},>=latex] 
\node[point] (1a) at (0,0) {};
\node[epoint] (1b) at (0,1) {};
\node[epoint] (1c) at (0,2) {};
\node[point] (2a) at (1,0) {};
\node[point] (2b) at (1,1) {};
\node[point] (2c) at (1,2) {};
\node[epoint] (3a) at (2,0) {};
\node[point] (3b) at (2,1) {};
\node[point] (3c) at (2,2) {};
\node[point] (4b) at (3,1) {};
\node[epoint] (4c) at (3,2) {};
  
\path[->,thick]
   (2a) edge[red,dashed] (1b)
   (2b) edge[red,dashed] (1c)
   (2c) edge (1a)
   (2a) edge[red,dashed] (3a)
   (2b) edge (3c)
   (2c) edge (3b)
   (4b) edge (3b)
   (4c) edge[red,dashed] (3c);
  \end{tikzpicture}
\qquad
\begin{tikzpicture}[point/.style={shape=circle,fill=black,scale=.5pt,outer sep=3pt},epoint/.style={shape=rectangle,fill=red,scale=.5pt,outer sep=3pt},>=latex] 
\node[point] (1a) at (0,0) {};
\node[epoint] (1b) at (0,1) {};
\node[epoint] (1c) at (0,2) {};
\node[point] (2a) at (1,0) {};
\node[point] (2b) at (1,1) {};
\node[point] (2c) at (1,2) {};
\node[epoint] (3a) at (2,0) {};
\node[point] (3b) at (2,1) {};
\node[point] (3c) at (2,2) {};
\node[point] (4b) at (3,1) {};
\node[epoint] (4c) at (3,2) {};
  
\path[->,thick]
   (2a) edge[red,dashed] (1b)
   (2b) edge[red,dashed] (1c)
   (2c) edge (1a)
   (2a) edge[red,dashed] (3a)
   (2b) edge (3b)
   (2c) edge (3c)
   (4b) edge (3c)
   (4c) edge[red,dashed] (3b);
  \end{tikzpicture}
\qquad
\begin{tikzpicture}[point/.style={shape=circle,fill=black,scale=.5pt,outer sep=3pt},epoint/.style={shape=rectangle,fill=red,scale=.5pt,outer sep=3pt},>=latex] 
\node[point] (1a) at (0,0) {};
\node[epoint] (1b) at (0,1) {};
\node[epoint] (1c) at (0,2) {};
\node[point] (2a) at (1,0) {};
\node[point] (2b) at (1,1) {};
\node[point] (2c) at (1,2) {};
\node[epoint] (3a) at (2,0) {};
\node[point] (3b) at (2,1) {};
\node[point] (3c) at (2,2) {};
\node[point] (4b) at (3,1) {};
\node[epoint] (4c) at (3,2) {};
  
\path[->,thick]
   (2a) edge[red,dashed] (1b)
   (2b) edge (1a)
   (2c) edge[red,dashed] (1c)
   (2a) edge (3c)
   (2b) edge[red,dashed] (3a)
   (2c) edge (3b)
   (4b) edge (3b)
   (4c) edge[red,dashed] (3c);

  \end{tikzpicture}\\
\vspace{1cm}
\begin{tikzpicture}[point/.style={shape=circle,fill=black,scale=.5pt,outer sep=3pt},epoint/.style={shape=rectangle,fill=red,scale=.5pt,outer sep=3pt},>=latex] 
\node[point] (1a) at (0,0) {};
\node[epoint] (1b) at (0,1) {};
\node[epoint] (1c) at (0,2) {};
\node[point] (2a) at (1,0) {};
\node[point] (2b) at (1,1) {};
\node[point] (2c) at (1,2) {};
\node[epoint] (3a) at (2,0) {};
\node[point] (3b) at (2,1) {};
\node[point] (3c) at (2,2) {};
\node[point] (4b) at (3,1) {};
\node[epoint] (4c) at (3,2) {};
  
\path[->,thick]
   (2a) edge[red,dashed] (1b)
   (2b) edge (1a)
   (2c) edge[red,dashed] (1c)
   (2a) edge[red,dashed] (3a)
   (2b) edge (3c)
   (2c) edge (3b)
   (4b) edge (3c)
   (4c) edge[red,dashed] (3b);

  \end{tikzpicture}
\qquad
\begin{tikzpicture}[point/.style={shape=circle,fill=black,scale=.5pt,outer sep=3pt},epoint/.style={shape=rectangle,fill=red,scale=.5pt,outer sep=3pt},>=latex] 
\node[point] (1a) at (0,0) {};
\node[epoint] (1b) at (0,1) {};
\node[epoint] (1c) at (0,2) {};
\node[point] (2a) at (1,0) {};
\node[point] (2b) at (1,1) {};
\node[point] (2c) at (1,2) {};
\node[epoint] (3a) at (2,0) {};
\node[point] (3b) at (2,1) {};
\node[point] (3c) at (2,2) {};
\node[point] (4b) at (3,1) {};
\node[epoint] (4c) at (3,2) {};
  
\path[->,thick]
   (2a) edge (1a)
   (2b) edge[red,dashed] (1b)
   (2c) edge[red,dashed] (1c)
   (2a) edge (3c)
   (2b) edge[red,dashed] (3a)
   (2c) edge (3b)
   (4b) edge (3c)
   (4c) edge[red,dashed] (3b);
  \end{tikzpicture}
\qquad
\begin{tikzpicture}[point/.style={shape=circle,fill=black,scale=.5pt,outer sep=3pt},epoint/.style={shape=rectangle,fill=red,scale=.5pt,outer sep=3pt},>=latex] 
\node[point] (1a) at (0,0) {};
\node[epoint] (1b) at (0,1) {};
\node[epoint] (1c) at (0,2) {};
\node[point] (2a) at (1,0) {};
\node[point] (2b) at (1,1) {};
\node[point] (2c) at (1,2) {};
\node[epoint] (3a) at (2,0) {};
\node[point] (3b) at (2,1) {};
\node[point] (3c) at (2,2) {};
\node[point] (4b) at (3,1) {};
\node[epoint] (4c) at (3,2) {};
  
\path[->,thick]
   (2a) edge[red,dashed] (1b)
   (2b) edge (1a)
   (2c) edge[red,dashed] (1c)
   (2a) edge (3b)
   (2b) edge[red,dashed] (3a)
   (2c) edge (3c)
   (4b) edge (3c)
   (4c) edge[red,dashed] (3b);
  \end{tikzpicture}
    \caption{3 crossing $K$-theoretic diagrams for $\cO_{\bw}$}
\label{fig:laces4}
\end{figure}

\begin{figure}
\begin{tikzpicture}[point/.style={shape=circle,fill=black,scale=.5pt,outer sep=3pt},epoint/.style={shape=rectangle,fill=red,scale=.5pt,outer sep=3pt},>=latex,scale=0.8]
\node[point] (1a) at (0,0) {};
\node[epoint] (1b) at (0,1) {};
\node[epoint] (1c) at (0,2) {};
\node[point] (2a) at (1,0) {};
\node[point] (2b) at (1,1) {};
\node[point] (2c) at (1,2) {};
\node[epoint] (3a) at (2,0) {};
\node[point] (3b) at (2,1) {};
\node[point] (3c) at (2,2) {};
\node[point] (4b) at (3,1) {};
\node[epoint] (4c) at (3,2) {};
  
\path[->,thick]
   (2a) edge[red,dashed] (1b)
   (2b) edge[red,dashed] (1c)
   (2c) edge (1a)
   (2a) edge (3c)
   (2b) edge[red,dashed] (3a)
   (2c) edge (3b)
   (4b) edge (3b)
   (4c) edge[red,dashed] (3c);
  \end{tikzpicture}
\quad
\begin{tikzpicture}[point/.style={shape=circle,fill=black,scale=.5pt,outer sep=3pt},epoint/.style={shape=rectangle,fill=red,scale=.5pt,outer sep=3pt},>=latex,scale=0.8]
\node[point] (1a) at (0,0) {};
\node[epoint] (1b) at (0,1) {};
\node[epoint] (1c) at (0,2) {};
\node[point] (2a) at (1,0) {};
\node[point] (2b) at (1,1) {};
\node[point] (2c) at (1,2) {};
\node[epoint] (3a) at (2,0) {};
\node[point] (3b) at (2,1) {};
\node[point] (3c) at (2,2) {};
\node[point] (4b) at (3,1) {};
\node[epoint] (4c) at (3,2) {};
  
\path[->,thick]
   (2a) edge[red,dashed] (1b)
   (2b) edge[red,dashed] (1c)
   (2c) edge (1a)
   (2a) edge[red,dashed] (3a)
   (2b) edge (3c)
   (2c) edge (3b)
   (4b) edge (3c)
   (4c) edge[red,dashed] (3b);
  \end{tikzpicture}
\quad
\begin{tikzpicture}[point/.style={shape=circle,fill=black,scale=.5pt,outer sep=3pt},epoint/.style={shape=rectangle,fill=red,scale=.5pt,outer sep=3pt},>=latex,scale=0.8]
\node[point] (1a) at (0,0) {};
\node[epoint] (1b) at (0,1) {};
\node[epoint] (1c) at (0,2) {};
\node[point] (2a) at (1,0) {};
\node[point] (2b) at (1,1) {};
\node[point] (2c) at (1,2) {};
\node[epoint] (3a) at (2,0) {};
\node[point] (3b) at (2,1) {};
\node[point] (3c) at (2,2) {};
\node[point] (4b) at (3,1) {};
\node[epoint] (4c) at (3,2) {};
  
\path[->,thick]
   (2a) edge[red,dashed] (1b)
   (2b) edge (1a)
   (2c) edge[red,dashed] (1c)
   (2a) edge (3c)
   (2b) edge[red,dashed] (3a)
   (2c) edge (3b)
   (4b) edge (3c)
   (4c) edge[red,dashed] (3b);

  \end{tikzpicture}
\quad
\begin{tikzpicture}[point/.style={shape=circle,fill=black,scale=.5pt,outer sep=3pt},epoint/.style={shape=rectangle,fill=red,scale=.5pt,outer sep=3pt},>=latex,scale=0.8]
\node[point] (1a) at (0,0) {};
\node[epoint] (1b) at (0,1) {};
\node[epoint] (1c) at (0,2) {};
\node[point] (2a) at (1,0) {};
\node[point] (2b) at (1,1) {};
\node[point] (2c) at (1,2) {};
\node[epoint] (3a) at (2,0) {};
\node[point] (3b) at (2,1) {};
\node[point] (3c) at (2,2) {};
\node[point] (4b) at (3,1) {};
\node[epoint] (4c) at (3,2) {};
  
\path[->,thick]
   (2a) edge[red,dashed] (1b)
   (2b) edge[red,dashed] (1c)
   (2c) edge (1a)
   (2a) edge (3c)
   (2b) edge[red,dashed] (3a)
   (2c) edge (3b)
   (4b) edge (3c)
   (4c) edge[red,dashed] (3b);
  \end{tikzpicture}
    \caption{4 or 5 crossing  $K$-theoretic diagrams for $\cO_{\bw}$}
\label{fig:laces5}
\end{figure}
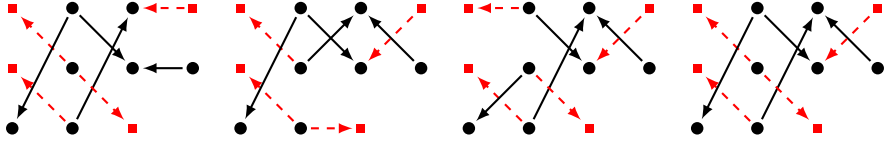

There are also \emph{$K$-theoretic transformations} of lacing diagrams
\begin{equation}
\vcenter{\hbox{  
\begin{tikzpicture}[point/.style={shape=circle,fill=black,scale=.5pt,outer sep=3pt},>=latex]
\node[point] (1a) at (0,0) {};
\node[point] (1b) at (0,1) {};
\node[point] (2a) at (1,0) {};
\node[point] (2b) at (1,1) {};
\node[point] (3a) at (2,0) {};
\node[point] (3b) at (2,1) {};

\path[-,thick]
   (1a) edge (2b)
   (1b) edge (2a)
   (2a) edge (3a)
   (2b) edge (3b);
\end{tikzpicture}}}
\longleftrightarrow
\vcenter{\hbox{  
\begin{tikzpicture}[point/.style={shape=circle,fill=black,scale=.5pt,outer sep=3pt},>=latex]
\node[point] (1a) at (0,0) {};
\node[point] (1b) at (0,1) {};
\node[point] (2a) at (1,0) {};
\node[point] (2b) at (1,1) {};
\node[point] (3a) at (2,0) {};
\node[point] (3b) at (2,1) {};

\path[-,thick]
   (1a) edge (2b)
   (1b) edge (2a)
   (2a) edge (3b)
   (2b) edge (3a);
\end{tikzpicture}}}
\longleftrightarrow
\vcenter{\hbox{  
\begin{tikzpicture}[point/.style={shape=circle,fill=black,scale=.5pt,outer sep=3pt},>=latex]
\node[point] (1a) at (0,0) {};
\node[point] (1b) at (0,1) {};
\node[point] (2a) at (1,0) {};
\node[point] (2b) at (1,1) {};
\node[point] (3a) at (2,0) {};
\node[point] (3b) at (2,1) {};

\path[-,thick]
   (1a) edge (2a)
   (1b) edge (2b)
   (2a) edge (3b);
   (2b) edge (3a)
\end{tikzpicture}}}
\end{equation}
with the same condition on the dots, and in addition the two middle dots should be consecutive in their column.  Suppose $\bw$ is a minimal lacing diagram.  A lacing diagram $\bw'$ is said to be a \emph{$K$-theoretic lacing diagram} for the orbit $\cO_\bw$ if $\bw'$ can be obtained from $\bw$ by a sequence of $K$-theoretic transformations of diagrams.  The $K$-theoretic lacing diagrams for the running example are shown in Figures \ref{fig:laces4} and \ref{fig:laces5}.  Note that last lacing diagram in Figure \ref{fig:laces2} lies in $\cO_{\bw}$ but is not a $K$-theoretic lacing diagram for that orbit.

We see in the above example that it is not necessarily true that a $K$-theoretic lacing diagram for $\cO_\bw$ will lie in the same orbit as $\bw$.  We will see below that non-minimal $K$-theoretic diagrams for $\cOc{\bw}$ correspond to higher degree terms in the $K$-polynomial of $\cOc{\bw}$, arising from an inclusion-exclusion process on $K$-classes of irreducible components of a partial Gr\"obner degeneration of $\cOc{\bw}$.  A purely representation theoretic meaning of $K$-theoretic diagrams is not clear.

\section{K-polynomials of quiver orbit closures}
We start this section by recalling the definition of $K$-polynomials in the case of closed subvarieties of an affine space. Fix and action of the algebraic torus $T=(\kk^\times)^n$ on affine space $\mathbb{A}^N$ and let $X \subseteq \mathbb{A}^N$ be a $T$-stable closed subvariety.  This induces a $\mathbb{Z}^n$-grading on the coordinate ring
\begin{equation}
\kk[X] = \bigoplus_{\mathbf{e} \in \mathbb{Z}^n} \kk[X]_\mathbf{e}, 
\end{equation}
where the homogeneous pieces are the weight spaces (see, for example, \cite[\S2.6]{KKR} for a detailed account).  We assume for the remainder of the discussion that this grading is \emph{positive}, meaning that the degree $\mathbf{0}$ piece consists only of the coefficient field $\kk$ (which will always be the case in our application to quiver representations below). This grading gives rise to a Hilbert series, defined below, which lives in the additive group $\prod_{\mathbf{e} \in \mathbb{Z}^n} \mathbb{Z} \mathbf{a}^\mathbf{e}$, where $\mathbf{a}=(a_1, \dotsc, a_n)$ is an alphabet,  $\mathbf{e} = (e_1, \dotsc, e_n) \in \mathbb{Z}^n$, and $\mathbf{a}^{\mathbf{e}} = a_1^{e_1} \cdots a_n^{e_n}$ is a monomial. 

The (multigraded) \emph{Hilbert series} of $X$ is defined as
\begin{equation}
H(X; \mathbf{a}) = \sum_{\mathbf{e} \in \mathbb{Z}^n} \left(\dim_\kk \kk[X]_\mathbf{e} \right)\mathbf{a}^\mathbf{e}.
\end{equation}
This important invariant has a drawback which we want to circumvent: suppose that $I \subseteq \kk[\mathbb{A}^N]$ is the defining ideal of $X$, and that we extend scalars to $I \otimes_\kk \kk[\mathbb{A}^{N'}] \subseteq \kk[\mathbb{A}^{N+N'}]$ to work in a larger ambient space (e.g., take the ideal generated by minors in a matrix of variables, then consider the ideal generated by the same minors in a larger matrix).  This will change the Hilbert series.
On the other hand, it can be shown that the ratio
\begin{equation}
\mathcal{K}(X; \mathbf{a}):= \frac{H(X; \mathbf{a})}{H(\mathbb{A}^N; \mathbf{a})}
\end{equation}
remains unchanged by extending scalars as above.  This expression is actually a Laurent polynomial in the variables $\mathbf{a}$, known as the \emph{$K$-polynomial of $\cOc{M}$}.  

For the remainder of this article, the varieties we work with will always come equipped with a natural torus action (which the reader will be reminded of), and thus it is unambiguous to omit the variable set $\mathbf{a}$ from the notation.  We will follow this practice in order to significantly simplify the notation.


For an arbitrary quiver $Q$ and dimension vector $\bd$, we consider the (maximal) torus $T \subseteq GL(\bd)$ consisting of all collections of diagonal matrices.  Then $\rep_Q(\bd)$ and all $\cOc{M}$ in it inherit an action of $T$ from this inclusion, and thus the coordinate ring $\kk[\cOc{M}]$ has a natural $\mathbb{Z}^d$-grading, where $d=\sum_{z \in Q_0} \bd(z)$.  Following the convention of the last paragraph, we may simply denote the $K$-polynomial of an orbit closure with respect to this torus action by $\mathcal{K}(\cOc{M})$.

In the running example, we make the identification of $\rep_Q(\bd)$ with the product of matrix spaces whose general element is shown in \eqref{eq:rep}.  The row and column labels illustrate the alphabets associated to the tori acting by row and column scaling.  Thus, the degree of a coordinate function picking out a matrix entry is its row label minus its column label; for example the degree of the coordinate function picking out the entry $b_3$ is $u_2 - s_1$.
\begin{equation}\label{eq:rep}
\rep_Q(\bd) = \left\{\left(
\begin{blockarray}{cc}
 & v_1\\
 \begin{block}{c[c]}
u_1 & c_1 \\
u_2 & c_2 \\
u_3 & c_3 \\
\end{block}
\end{blockarray},
\begin{blockarray}{ccc}
 & s_1 & s_2 \\
 \begin{block}{c[cc]}
u_1 & b_1 & b_2 \\
u_2 & b_3 & b_4 \\
u_3 & b_5 & b_6 \\
\end{block}
\end{blockarray},
\begin{blockarray}{ccc}
 & s_1 & s_2 \\
 \begin{block}{c[cc]}
t_1 & a_1 & a_2 \\
\end{block}
\end{blockarray}
\right) \right\}.
\end{equation}
The orbit closure of the lacing diagrams of our running example (Figure \ref{fig:laces}) is defined by the rank condition
\begin{equation}\label{eq:ranks}
\rank \begin{bmatrix}
0 & a_1 & a_2 \\
c_1 & b_1 & b_2 \\
c_2 & b_3 & b_4 \\
c_3 & b_5 & b_6 \\
\end{bmatrix} \leq 2,
\end{equation}
so that the ideal of this orbit closure is generated by the $3\times3$-minors of the matrix in \eqref{eq:ranks}.

\begin{remark}
The invariant $\mathcal{K}(\cOc{M})$ has several equivalent formulations, whose relations are carefully explained in \cite[\S\S 3,4]{MR2492443}.  For example, it represents an element of the ring of virtual rational representations of $GL(\bd)$.  This ring in turn can be identified with the Grothendieck group of the category of $GL(\bd)$-equivariant coherent sheaves on $\rep_Q(\bd)$, known as the $GL(\bd)$-equivariant $K$-homology of $\rep_Q(\bd)$, which is isomorphic to the $GL(\bd)$-equivariant $K$-cohomology ring of $\rep_Q(\bd)$ since this variety is nonsingular.
\end{remark}

The building blocks of our formulas for $K$-polynomials of quiver orbit closures are the double Grothendieck polynomials $\mathfrak{G}_w(\mathbf{a}; \mathbf{b})$ of Lascoux and Sch\"utzenberger.
Instead of the original recursive definition (see for example \cite{LS,fultonlascoux}), we will introduce them below as $K$-polynomials of certain closed subvarieties of matrix spaces, since this is how they naturally arise in our work and significantly simplifies the overview.
The connection between this formulation and the original definition is explained in \cite[\S2.7]{KKR}, following from results in \cite{MR1932326} or \cite{KM05}.

Given a space of matrices $X=\Mat(m,n)$, consider the natural multiplication action of $B_- \times B_+$ on $X$ where $B_- \subseteq GL(m)$ denotes the group of invertible lower triangular matrices, and $B_+ \subseteq GL(n)$  the group of invertible upper triangular matrices.  Each $B_- \times B_+$-orbit on $X$ has a unique partial permutation matrix; denote by $X_w$ the closure of the orbit containing the partial permutation matrix $w$, which is called a \emph{matrix Schubert variety}.  The equations defining $X_w$ as a closed subscheme are collections of minors corresponding to imposing upper bounds on the ranks of all northwest justified submatrices in the space $X$.  
A more detailed introduction to matrix Schubert varieties and their properties can be found in \cite[Ch. 15]{MS}.  

Each matrix Schubert variety $X_w \subseteq X$ carries an action of the subgroup of diagonal matrices
\begin{equation}
(\kk^\times)^m \times (\kk^\times)^n \subseteq B_- \times B_+
\end{equation}
where the factor $(\kk^\times)^m$ acts by scaling rows, and the factor $(\kk^\times)^n$ acts by scaling columns.  The coordinate ring $\kk[X_w]$ thus inherits a $\mathbb{Z}^{m+n}$ grading, and we identify this grading group with the free abelian group on the alphabet $(\mathbf{a}, \mathbf{b}) := (a_1, \dotsc, a_m, b_1, \dotsc, b_n)$.  Then we will take
\begin{equation}
\mathfrak{G}_w(\mathbf{a}; \mathbf{b}) = \mathcal{K}(X_w; \mathbf{a}, \mathbf{b})
\end{equation}
to be the \emph{double Grothendieck polynomial} indexed by $w$.  As mentioned above, we will omit the variables and simply denote it by $\mathfrak{G}_w$ below.  In addition to the original definition mentioned above, there are many other combinatorial formulas for Grothendieck polynomials \cite{MR2307216,MR1680646,MR1763950,BRspec,MR1946917}.

To connect with our running example, consider the leftmost lacing diagram in Figure \ref{fig:laces3}; written in matrix form we find the middle matrix to be 
\begin{equation}
w_b =\begin{blockarray}{ccc}
 & s_1 & s_2 \\
 \begin{block}{c[cc]}
u_1 & 0 & 1 \\
u_2 & 0 & 0 \\
u_3 & 1 & 0 \\
\end{block}
\end{blockarray}.
\end{equation}
The corresponding Grothendieck polynomial can be calculated (for example, using Lascoux-Sch\"utzenberger's recursive definition or one of the other combinatorial formulas cited above) to be
\begin{equation}
\mathfrak{G}_{w_b} = \left(1 - \frac{u_1}{s_1}\right)\left(1-\frac{u_2}{s_1}\right).
\end{equation}

Restricting our attention to type $A$ quivers now, we need to introduce some ``opposite'' versions of the above concepts to deal with arbitrary orientation.
For each partial permutation matrix $w$, the \emph{opposite matrix Schubert variety} $X^w \subseteq X$ is the closure of the $B_+ \times B_-$-orbit containing $w$, where $B_+ \subseteq GL(m)$ denotes the group of invertible upper triangular matrices, and $B_- \subseteq GL(n)$  the group of invertible lower triangular matrices.
Likewise, we get an \emph{opposite Grothendieck polynomial}
\begin{equation}
\mathfrak{G}^w(\mathbf{a};\mathbf{b}) = \mathcal{K}(X^w; \mathbf{a}, \mathbf{b}),
\end{equation}
denoted $\mathfrak{G}^w$ for short.  If $w$ is a permutation matrix, then $\mathfrak{G}^w$ is just the standard double Grothendieck polynomial for the $180^\circ$ rotation of $w$, with the orders of the individual input alphabets $\mathbf{a},\mathbf{b}$ reversed.

Then a lacing diagram $\bw=(w_a)_{a \in Q_1}$ determines a product of matrix Schubert varieties and opposite matrix Schubert varieties
\begin{equation}\label{eq:Xbw}
X_\bw = \prod_{\xrightarrow{a}\in Q_1} X_{w_a} \times \prod_{\xleftarrow{a}\in Q_1} X^{w_a} \subseteq \rep_Q(\bd) 
\end{equation}
where the first product is over rightward arrows of $Q$ and the second product over leftward arrows of $Q$. 
This is a $T$-stable subvariety and the $K$-polynomial of $X_\bw$ is simply the following product of Grothendieck polynomials (and their opposites):
\begin{equation}
\mathfrak{G}_\bw = \left(\prod_{\xrightarrow{a}\in Q_1} \mathfrak{G}_{w_a} \right) \left( \prod_{\xleftarrow{a}\in Q_1} \mathfrak{G}^{w_a} \right).
\end{equation}
These are the building blocks of our main result below.  For example, consider the unique 5-crossing $K$-theoretic lace diagram for the running example, which is the rightmost entry of Figure \ref{fig:laces5}.  Its matrix representation is
\begin{equation}
\bw =(w_a, w_b, w_c) = \left( 
\begin{blockarray}{cc}
 & v_1\\
 \begin{block}{c[c]}
u_1 & 1 \\
u_2 & 0 \\
u_3 & 0 \\
\end{block}
\end{blockarray},
\begin{blockarray}{ccc}
 & s_1 & s_2 \\
 \begin{block}{c[cc]}
u_1 & 0 & 1 \\
u_2 & 0 & 0 \\
u_3 & 1 & 0 \\
\end{block}
\end{blockarray},
\begin{blockarray}{ccc}
 & s_1 & s_2 \\
 \begin{block}{c[cc]}
t_1 & 1 & 0 \\
\end{block}
\end{blockarray}
\right)
\end{equation}
and the corresponding product of Grothendieck polynomials is
\begin{equation}\label{eq:lacepoly}
\mathfrak{G}_\bw = \mathfrak{G}^{w_c}\mathfrak{G}_{w_b}\mathfrak{G}^{w_a} = \left(1 - \frac{u_2}{v_1}\right)\left(1-\frac{u_3}{v_1}\right) \left(1 - \frac{u_1}{s_1}\right)\left(1-\frac{u_2}{s_1}\right)\left( 1- \frac{t_1}{s_2}\right).
\end{equation}

\section{The component formula}
At this point, we have enough background in place to state our capstone formula from \cite{KKR} for $K$-polynomials of type $A$ orbit closures.  This formula is called the ``component formula'' since it generalizes a formula of the same name from \cite{KMS}, where it was proven in the case of equioriented type $A$ quivers.  We remark that the name comes from the geometry technique of the proof outlined below; the name ``lace formula'' would suit it equally well.  A detailed account of the motivations for this formula and its relation to existing literature can be found in Sections 1.1 and 1.3 of \cite{KKR}, respectively.  Let us at least mention here though that it was conjectured by Buch and Rim\'anyi in \cite{MR2306279}, where they proved the cohomological version using the \emph{interpolation method} of Feh\'er and Rim\'anyi \cite{MR2087805,BFR}.

\begin{theorem}[Theorems 4.37 and 5.20 of \cite{KKR}]\label{thm:main}
Let $Q$ be an arbitrary quiver of Dynkin type $A$, let $\bd$ be a dimension vector for $Q$, and $\cOc{M} \subseteq \rep_Q(\bd)$ an orbit closure.  Then the $K$-polynomial of $\cOc{M}$ is given by the formula
\begin{equation}\label{eq:componentformula}
\mathcal{K}(\cOc{M}) = \sum_{\bw} (-1)^{|\bw| - {\rm codim}\cOc{M}} \mathfrak{G}_\bw
\end{equation}
where the sum is over $K$-theoretic lacing diagrams for $\cOc{M}$.
\end{theorem}

So in our running example we would sum 15 polynomials of the form \eqref{eq:lacepoly} in 12 variables, indexed by the 15 $K$-theoretic lacing diagrams shown in the figures of Section \ref{sec:lacing}.

We will now give a bird's eye view of the proof of Theorem \ref{thm:main}.
The first main idea is to reduce the problem to type $A$ quivers of a specific orientation, namely the bipartite (i.e., sink-source) orientation.  Given $Q$ of type $A$, one simply inserts a ``backwards'' arrow in the middle of each length two path to get an associated bipartite type $A$ quiver $\tilde{Q}$.  Then there is a dimension vector $\widetilde{\bd}$ for $\widetilde{Q}$ such that the equivariant geometry of $\rep_Q(\bd)$ can be relatively easily reduced to that of $\rep_{\widetilde{Q}}(\widetilde{\bd})$.  It was pointed out to us by Jorge Vit\'oria that the path algebras of $Q$ and $\widetilde{Q}$ are related by a universal localization in the sense of \cite[\S4]{MR800853}.  In contrast, several prominent results preceding ours on the geometry of orbit closures for arbitrary type $A$ quivers \cite{BZ-typeA,MR3008913} were reduced to the equioriented case using Zwara's work on Hom-controlled functors \cite{MR1888422} and sophisticated representation theoretical arguments.  The results in that case had been proven in work of Lakshmibai and Magyar \cite{LMdegen} using the connections with Schubert varieties mentioned below.

Reducing to the case of bipartite $Q$ allows us to make liberal use of the bipartite Zelevinsky map from the author's previous work with Rajchgot \cite{KR} throughout all parts of the proof.  This map embeds $\rep_Q(\bd)$ in a partial flag variety $GL(d)/P$, such that orbit closures in $\rep_Q(\bd)$ are identified with open subschemes of Schubert varieties.  This map generalizes results of Zelevinsky \cite{Zgradednilp} and Lakshmibai-Magyar \cite{LMdegen} on the case of equioriented type $A$ quivers, and in particular yields defining equations of the prime ideals of the orbit closures in that case.  While passing back and forth between $\rep_Q(\bd)$ and $GL(d)/P$ via this map is ubiquitous throughout the proof, it is primarily for technical purposes, and thus will not be discussed further.

The second main idea can be explained in a purely geometric way, though a significant amount of combinatorics is necessary for the proof. We construct a ``simultaneous flat degeneration" of the action of $GL(\bd)$ on $\cOc{M}$ to an action of $B_+ \times_T B_-$ on a scheme $\cOc{M}(0)$.
Here, $B_+ \times_T B_-$ denotes a certain subgroup of $GL(\bd) \times GL(\bd)$ consisting of various upper and lower triangular matrices, determined by the orientation of $Q$, and the degeneration $\cOc{M}(0)$ is a union of a certain (potentially nonreduced) subschemes whose underlying varieties are various $X_\bw \subseteq \rep_Q(\bd)$, as defined in \eqref{eq:Xbw} above.  The key point of this stage of our proof is that
\begin{equation}
\cOc{M}(0) = \bigcup_{\bw\ \text{minimal}} X_\bw
\end{equation}
(i.e., the degeneration is reduced and its irreducible components are precisely the $B_+ \times_T B_-$-orbit closures of the minimal lace diagrams for $M$).  Because of flatness, the $K$-polynomial of $\cOc{M}$ is equal to the $K$-polynomial of $\cOc{M}(0)$. 

We remark that this is already enough to prove the cohomological component formula of \cite{MR2306279}, since this invariant only requires knowledge of the irreducible components of a scheme.  On the other hand, the $K$-polynomial contains deeper information depending on how the irreducible components intersect.  A simple but instructive comparison can be found in \cite[Example 1]{MR2143073}.  Thus, roughly speaking, it remains to understand the configuration of the irreducible components of the degeneration with respect to one another.  

\begin{figure}
\begin{tikzpicture}[c/.style={shape=circle,draw}]
\node[c] (a) at (1, 7) {$+1$};
\node[c] (b) at (3, 7) {$+1$};
\node[c] (c) at (5, 7) {$+1$};
\node[c] (d) at (7, 7) {$+1$};
\node[c] (e) at (9, 7) {$+1$};
\node[c] (ab) at (0, 5) {$-1$};
\node[c] (ad) at (2, 5) {$-1$};
\node[c] (bc) at (4, 5) {$-1$};
\node[c] (bd) at (6, 5) {$-1$};
\node[c] (ce) at (8, 5) {$-1$};
\node[c] (de) at (10, 5) {$-1$};
\node[c] (z) at (3, 3) {$+1$};
\node[c] (x) at (5, 3) {$+1$};
\node[c] (y) at (7, 3) {$+1$};
\node[c] (w) at (5, 1) {$-1$};
\draw[-]
 (a) edge (ab)
 (a) edge (ad)
 (b) edge (ab)
 (b) edge (bc)
 (b) edge (bd)
 (c) edge (bc)
 (c) edge (ce)
 (d) edge (ad)
 (d) edge (bd)
 (d) edge (de)
 (e) edge (ce)
 (e) edge (de);
\draw[-]
 (ab) edge (x)
 (ab) edge (z)
 (ad) edge (x)
 (bc) edge (y)
 (bc) edge (z)
 (bd) edge (x)
 (bd) edge (y)
 (ce) edge (y) 
 (de) edge (y);
\draw[-]
(x) edge (w)
(y) edge (w)
(z) edge (w);
\end{tikzpicture}
    \caption{The Hasse diagram of the poset of $K$-theoretic diagrams for $\cO_{\bw}$ in the running example, with each element labeled by the value the M\"obius function takes there}
\label{fig:hasse}
\end{figure}
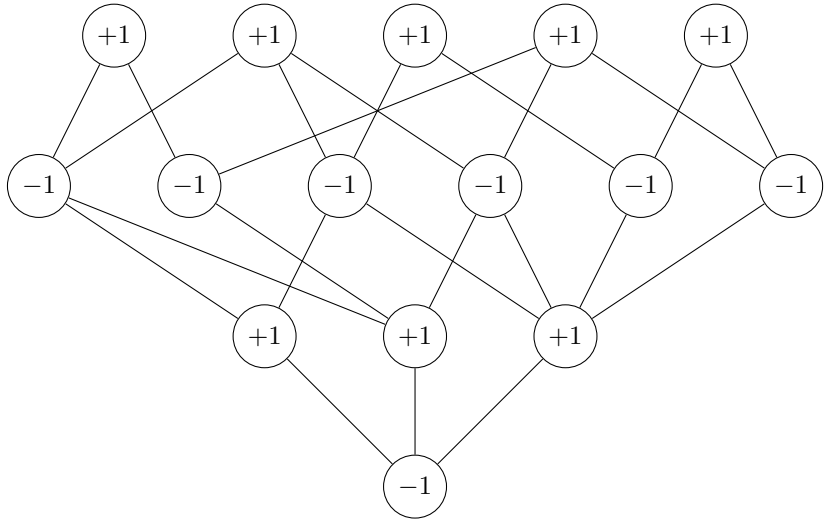

The third and final main idea is a computation of the M\"obius functions of certain posets in terms of lacing diagram combinatorics.  To be more precise, a theorem of Knutson \cite{Knutson:2009aa} allows us to compute the $K$-polynomial of $\cOc{M}(0)$ by computing the values of the M\"obius function of the poset $\mathcal{M}$
of $B_+ \times_T B_-$-orbit closures contained in $\cOc{M}(0)$, with respect to containment order.  The goal is to show that the M\"obius function is nonzero on an orbit closure $X_\bw$ precisely when $\bw$ is a $K$-theoretic lacing diagram for $\cO_M$, and in this case the value is $\pm 1$, alternating with the number of crossings of $\bw$.  In Figure \ref{fig:hasse} we see the Hasse diagram of the poset of $K$-theoretic lacing diagrams for $\cO_M$ in our running example, ordered by containment of their corresponding $B_+ \times_T B_-$-orbit closures, with each element labeled by the value of the M\"obius function.  Starting at the top left and reading each row to the right, then proceeding down to the next row, the ordering is the same as the order in which they are displayed in the figures of Section \ref{sec:lacing}.

This M\"obius function is difficult to understand directly, so we embed the incidence algebra of $\mathcal{M}$ in the incidence algebra of
the \emph{pipe complex} $\Delta$ studied by Woo and Yong \cite{wooyong2} (which is an example of the subword complexes studied by Knutson and Miller \cite{MR2047852}).  It is easier to compute the M\"obius function of $\Delta$ because this poset is a simplicial complex which is homeomorphic to a ball.  Then the primary challenge of this step is 
to relate the combinatorics of the pipe complex to those of the $K$-theoretic lacing diagrams in order to obtain the desired results on the M\"obius function of $\mathcal{M}$, thus completing the proof.

\section{Open problems}

Finally, we collect some ideas for future directions to be pursued.  The reader is also referred to Zwara's survey \cite{Zwarasurvey} for a wealth of interesting problems and questions about the geometry of orbit closures in $\rep_Q(\bd)$ in more generality.

\subsubsection*{Formulas for $K$-polynomials}
For an arbitrary quiver $Q$ without oriented cycles, Buch has given a general shape of formulas for $K$-polynomials of $GL(\bd)$-invariant closed subvarieties $\Omega \subseteq \rep_Q(\bd)$ (e.g., orbit closures) in \cite{MR2492443}. 
He showed that each such $K$-polynomial can be written as a certain sum of products of \emph{stable double Grothendieck polynomials} indexed by partitions.  His expression has a uniqueness property that leads to well-defined \emph{quiver coefficients} $c_\mu(\Omega) \in \mathbb{Z}$ indexed by sequences of partitions $\mu =(\mu_i)_{i \in Q_0}$.  His conjecture on the properties of these numbers is only fully proven for equioriented type $Q$ quivers, and quivers of type $A_3$.  Given the formal similarity of his formula to our component formula, it would be interesting to better understand their relation.

\begin{problem}
Use the component formula \eqref{eq:componentformula} to prove Buch's Conjecture 1.1 of \cite{MR2492443} for all type $A$ quivers.
\end{problem}


\subsubsection*{Embeddings of representation varieties}

Bobi\'nski and Zwara have shown in \cite{MR1967381} that the smooth equivalence classes of singularities which appear in type $A$ quiver orbit closures (varying over all type $A$ quivers $Q$ and $\bd$) are precisely those which appear in orbit closures in flag varieties $GL(n)/B$ (Schubert varieties).
It has been reported by Andr\'as L\H{o}rincz that his work on $b$-functions of quiver semi-invariants \cite{Lorincz1,Lorincz2} implies that there exist singularities in type $D$ quiver orbit closures which are not smoothly equivalent to any singularity of any type $A$ quiver orbit closure.  Thus, we should look beyond the realm of flag varieties $GL(n)/B$ if we want to generalize the type $A$ Zelevinsky map and relate more general Dynkin quiver orbit closures to the wider algebraic geometry literature.  To retain as many of the desirable properties of the type $A$ Zelevinsky map as possible, \emph{spherical varieties} make reasonable candidates for target varieties; that is, varieties $X$ such that:
\begin{enumerate}
\item $X$ is normal (or even smooth);
\item $X$ comes equipped with the action of a connected, reductive algebraic group $G$;
\item $X$ has finitely many $B$-orbits, where $B\subset G$ is a Borel subgroup.
\end{enumerate}
These are well studied so that we may hope to transport their properties to $\rep_Q(\bd)$ (see the survey \cite{MR3177371}).

\begin{problem}
Let $Q$ be a quiver of Dynkin type $D$ or $E$.  Find a collection of varieties $X(\bd)$ along with maps $\zeta_\bd \colon \rep_Q(\bd) \to X(\bd)$ such that:
\begin{enumerate}
\item each $X(\bd)$ is spherical with respect to action of a reductive group $G(\bd)$, say with Borel subgroup $B(\bd)$;
\item each $\zeta_\bd$ is equivariant with respect to an embedding of algebraic groups $GL(\bd) \to G(\bd)$;
\item the restriction of $\zeta_\bd$ to an orbit closure gives an isomorphism between $\cOc{M}$ and the intersection of a $B(\bd)$-orbit closure with the image of $\zeta_\bd$, and these intersections differ from $B(\bd)$-orbit closures in $X(\bd)$ by a smooth factor.
\end{enumerate}
\end{problem}

Joint work in progress with Rajchgot, inspired by results of Bobi\'nski-Zwara and Brion \cite{MR1967381,MR2011763}, proposes a solution to this problem for $Q$ bipartite of type $D$ and each $X(\bd)$ a \emph{double Grassmannian} $Gr_a(V) \times Gr_b (V)$ with $G(\bd) = GL(V)$ acting diagonally ($a, b, V$ all depend on $\bd$).

\subsubsection*{Ideals of orbits closures}
The proofs in our paper \cite{KKR} produce a \emph{Gr\"obner basis} for the prime ideal of $\cOc{M}$.  Gr\"obner bases could be useful for attacking certain open problems, such as showing that a proposed generating set for the ideal of $\cOc{M}$ does indeed define a prime ideal, or for combinatorial cohomology formulas for $\cOc{M}$.
\begin{problem}
Let $Q$ be a Dynkin quiver and $\cOc{M}$ an orbit closure in some $\rep_Q(\bd)$.  Determine a Gr\"obner basis for the prime ideal defining $\cOc{M}$.
\end{problem}

\bibliographystyle{amsalpha}
\bibliography{ICRA}

\end{document}